\documentclass[12pt]{article}
\usepackage{graphicx}
\usepackage{amssymb,amsmath,amsfonts,amsthm}
\newtheorem{thm}{Theorem}[section]
\newtheorem{lem}[thm]{Lemma}        
\newtheorem{cor}[thm]{Corollary}

\newtheorem{conj}[thm]{Conjecture}
\newtheorem{obs}{Observation}

\newcommand{\qitem}[1]{\noindent\leavevmode\hangindent5mm%
       \noindent\hbox to5mm{#1\hss}\ignorespaces}
\newcommand{\qitemm}[1]{\noindent\leavevmode\hangindent10mm%
       \noindent\hbox to10mm{#1\hss}\ignorespaces}
\newcommand{\eop}{\qquad\hspace*{\fill} $\blacksquare$}
\newcommand{\openeop}{\qquad\hspace*{\fill} $\square$}
\newcommand{\defi}{\mathrm{def}}

\begin{document}

\title{\bf Degree Sequences and the Existence of $k$-Factors}

\author{{\sc  D. Bauer}\\
  {\small \sl Department of Mathematical Sciences} \\
  {\small \sl Stevens Institute of Technology } \\
  {\small \sl Hoboken, NJ 07030, U.S.A.}
  \and
  {\sc H.J. Broersma\,\footnote{\,Current address: Department of Applied
      Mathematics, University of Twente, P.O. Box 217, 7500 AE Enschede,
      The Netherlands}}\\
  {\small \sl Department of Computer Science}\\
  {\small \sl Durham University}\\
  {\small \sl South Road, Durham DH1 3LE, U.K.}
  \and
  {\sc J. van den Heuvel}\\
  {\small \sl Department of Mathematics}\\
  {\small \sl London School of Economics}\\
  {\small \sl Houghton Street, London WC2A 2AE, U.K.}
  \and
  {\sc N. Kahl}\\
  {\small \sl Department of Mathematics } \\
  {\small \sl and Computer Science}\\
  {\small \sl Seton Hall University}\\
  {\small \sl South Orange, NJ  07079, U.S.A.}
  \and
  {\sc E. Schmeichel}\\
  {\small \sl Department of Mathematics}\\
  {\small \sl San Jos\'e State University}\\
  {\small \sl San Jos\'e, CA  95192, U.S.A.}
}

\date{}
\maketitle

\begin{abstract}
\noindent
We consider sufficient conditions for a degree sequence $\pi$ to be
forcibly \mbox{$k$-factor} graphical. We note that previous work on degrees
and factors has focused primarily on finding conditions for a degree
sequence to be potentially $k$-factor graphical.

\vspace{1.5mm}\noindent
We first give a theorem for $\pi$ to be forcibly 1-factor graphical and,
more generally, forcibly graphical with deficiency at most $\beta \ge 0$.
These theorems are equal in strength to Chv\'atal's well-known
hamiltonian theorem, i.e., the best monotone degree condition for
hamiltonicity. We then give an equally strong theorem for $\pi$ to be
forcibly 2-factor graphical. Unfortunately, the number of nonredundant
conditions that must be checked increases significantly in moving from
$k=1$ to $k=2$, and we conjecture that the number of nonredundant
conditions in a best monotone theorem for a $k$-factor will increase
superpolynomially in $k$.

\vspace{1.5mm}\noindent
This suggests the desirability of finding a theorem for $\pi$ to be
forcibly \mbox{$k$-factor} graphical whose algorithmic complexity grows
more slowly. In the final section, we present such a theorem for any $k
\geq 2$, based on Tutte's well-known factor theorem. While this theorem is
not best monotone, we show that it is nevertheless tight in a precise way,
and give examples illustrating this tightness.


\end{abstract}

\section{Introduction}

We consider only undirected graphs without loops or multiple edges. Our
terminology and notation will be standard except as indicated, and a good
reference for any undefined terms or notation is \cite{CL}.

A \textit{degree sequence} of a graph on $n$ vertices is any sequence
$\pi=(d_1, d_2, \dots, d_n)$ consisting of the vertex degrees of the graph.
In contrast to \cite{CL}, we will usually assume the sequence is in
nondecreasing order. We generally use the standard abbreviated notation for
degree sequences, e.g., $(4,4,4,4,4,5,5)$ will be denoted $4^5 5^2$. A
sequence of integers $\pi = (d_1, d_2, \dots, d_n)$ is called
\textit{graphical} if there exists a graph~$G$ having $\pi$ as one of its
degree sequences, in which case we call $G$ a \textit{realization}
of~$\pi$. If $\pi = (d_1, \dots, d_n)$ and $\pi' = (d_1', \dots, d_n')$ are
two integer sequences, we say $\pi'$ \textit{majorizes}~$\pi$, denoted
$\pi' \geq \pi$, if $d_j' \geq d_j$ for $1 \leq j \leq n$. If $P$ is a
graphical property (e.g., $k$-connected, hamiltonian), we call a graphical
degree sequence \textit{forcibly} (respectively, \textit{potentially})
\textit{$P$ graphical} if every (respectively, some) realization of $\pi$
has property $P$.

Historically, the degree sequence of a graph has been used to provide
sufficient conditions for a graph to have a certain property, such as
$k$-connected or hamiltonian. Sufficient conditions for a degree sequence
to be forcibly hamiltonian were given by several authors, culminating in
the following theorem of Chv\'atal \cite{Chvatal72} in 1972.

\begin{thm}[\cite{Chvatal72}]\label{thm:chvatal}
  \; Let $\pi = (d_1 \leq \dots \leq d_n)$ be a graphical degree sequence,
  with $n \geq 3$. If $d_i \leq i < \frac12 n$ implies $d_{n-i} \geq n-i$,
  then $\pi$ is forcibly hamiltonian graphical.
\end{thm}

Unlike its predecessors, Chv\'atal's theorem has the property that if it
does not guarantee that a graphical degree sequence~$\pi$ is forcibly
hamiltonian graphical, then $\pi$ is majorized by some degree sequence
$\pi'$ which has a nonhamiltonian realization. As we'll see, this fact
implies that Chv\'atal's theorem is the strongest of an entire class of
theorems giving sufficient conditions for $\pi$ to be forcibly hamiltonian
graphical.

A \textit{factor} of a graph $G$ is a spanning subgraph of $G$. A
\textit{$k$-factor} of $G$ is a factor whose vertex degrees are identically
$k$. For a recent survey on graph factors, see \cite{Plummer07}. In the
present paper, we develop sufficient conditions for a degree sequence to be
forcibly $k$-factor graphical. We note that previous work relating degrees
and the existence of factors has focused primarily on sufficient conditions
for $\pi$ to be potentially \mbox{$k$-factor} graphical. The following
obvious necessary condition was conjectured to be sufficient by Rao and Rao
\cite{RaoRao72}, and this was later proved by Kundu \cite{Kundu73}.

\begin{thm}[\cite{Kundu73}]\label{kundu}
  \; The sequence $\pi = (d_1, d_2, \dots, d_n)$ is potentially $k$-factor
  graphical if and only if

  \qitemm{(1)}$(d_1, d_2, \dots, d_n)$ is graphical, and

  \qitemm{(2)}$(d_1-k, d_2-k, \dots, d_n-k)$ is graphical.
\end{thm}

Kleitman and Wang \cite{KW73} later gave a proof of Theorem \ref{kundu}
that yielded a polynomial algorithm constructing a realization $G$ of $\pi$
with a $k$-factor. Lov\'{a}sz \cite{Lovasz74} subsequently gave a very
short proof of Theorem \ref{kundu} for the special case $k=1$, and Chen
\cite{Chen88} produced a short proof for all $k \geq 1$.

In Section~\ref{s2}, we give a theorem for~$\pi$ to be forcibly graphical
with deficiency at most~$\beta$ (i.e., have a matching missing at
most~$\beta$ vertices), and show this theorem is strongest in the same
sense as Chv\'atal's hamiltonian degree theorem. The case $\beta=0$ gives
the strongest result for~$\pi$ to be forcibly 1-factor graphical. In
Section~\ref{s3}, we give the strongest theorem, in the same sense as
Chv\'atal, for~$\pi$ to be forcibly 2-factor graphical. But the increase in
the number of nonredundant conditions which must be checked as we move from
a 1-factor to a 2-factor is notable, and we conjecture the number of such
conditions in the best monotone theorem for $\pi$ to be forcibly $k$-factor
graphical increases superpolynomially in $k$. Thus it would be desirable to
find a theorem for $\pi$ to be forcibly $k$-factor graphical in which the
number of nonredundant conditions grows in a more reasonable way. In
Section \ref{s4}, we give such a theorem for $k \geq 2$, based on Tutte's
well-known factor theorem. While our theorem is not best monotone, it is
nevertheless tight in a precise way, and we provide examples to illustrate
this tightness.

We conclude this introduction with some concepts which are needed in the
sequel. Let $P$ denote a graph property (e.g., hamiltonian, contains a
\mbox{$k$-factor}, etc.) such that whenever a spanning subgraph of $G$ has
$P$, so does $G$. A function
$f : \{\text{Graphical Degree Sequences} \} \rightarrow \{0,1\}$ such that
$f(\pi)=1$ implies $\pi$ is forcibly $P$ graphical, and $f(\pi)=0$ implies
nothing in this regard, is called a \textit{forcibly $P$ function}. Such a
function is called \textit{monotone} if $\pi' \ge \pi$ and $f(\pi)=1$
implies $f(\pi')=1$, and \textit{weakly optimal} if $f(\pi)=0$ implies
there exists a graphical sequence $\pi' \ge \pi$ such that $\pi'$ has a
realization $G'$ without $P$. A forcibly $P$ function which is both
monotone and weakly optimal is the best monotone forcibly $P$ function, in
the following sense.

\begin{thm}
  \quad If $f$,$f_0$ are monotone, forcibly $P$ functions, and $f_0$ is
  weakly optimal, then $f_0(\pi) \ge f(\pi)$, for every graphical sequence
  $\pi$.
\end{thm}

\textbf{Proof:} Suppose to the contrary that for some graphical sequence
$\pi$ we have $1= f(\pi) > f_0(\pi) = 0$. Since $f_0$ is weakly optimal,
there exists a graphical sequence $\pi' \ge \pi$ such that $\pi'$ has a
realization $G'$ without $P$, and thus $f(\pi')=0$. But $\pi' \ge \pi$,
$f(\pi)=1$ and $f(\pi')=0$ imply $f$ cannot be monotone, a
contradiction.\eop

\bigskip
A theorem $T$ giving a sufficient condition for $\pi$ to be forcibly $P$
corresponds to the forcibly $P$ function $f_T$ given by: $f_T(\pi)=1$ if
and only if $T$ implies $\pi$ is forcibly~$P$. It is well-known that if $T$
is Theorem \ref{thm:chvatal} (Chv\'atal's theorem), then $f_T$ is both
monotone and weakly optimal, and thus the best monotone forcibly
hamiltonian function in the above sense. In the sequel, we will simplify
the formally correct `$f_T$ is monotone, etc.' to `$T$ is monotone, etc..'

\section{Best monotone condition for a 1-factor}\label{s2}

In this section we present best monotone conditions for a graph to have a
large matching. These results were first obtained by Las Vergnas
\cite{LV72}, and can also be obtained from results in Bondy and
Chv\'{a}tal~\cite{BC76}. For the convenience of the reader, we include the
statement of the results and short proofs below.

The \textit{deficiency} of $G$, denoted $\defi(G)$, is the number of
vertices unmatched under a maximum matching in $G$. In particular, $G$
contains a 1-factor if and only if $\defi(G)=0$.

We first give a best monotone condition for $\pi$ to be forcibly graphical
with deficiency at most $\beta$, for any $\beta \geq 0$.

\begin{thm}[\cite{BC76,LV72}]\label{thm:def}
  \; Let $G$ have degree sequence $\pi=(d_1 \leq \dots \leq d_n)$, and let
  $0 \leq \beta \leq n$ with $\beta \equiv n \pmod2$. If
  $$d_{i+1} \leq i-\beta < \tfrac12(n-\beta-1)\;\Longrightarrow\;
  d_{n+\beta-i} \geq n-i-1,$$
  then $\defi(G) \leq \beta$.
\end{thm}

The condition in Theorem \ref{thm:def} is clearly monotone. Furthermore, if
$\pi$ does not satisfy the condition for some $i \geq \beta$, then $\pi$ is
majorized by $\pi' = (i-\beta)^{i+1}\linebreak (n-i-2)^{n-2i+\beta-1}
(n-1)^{i-\beta}$. But $\pi'$ is realizable as
$K_{i-\beta}+(\overline{K_{i+1}} \cup K_{n-2i+\beta-1})$, which has
deficiency $\beta+2$. Thus Theorem \ref{thm:def} is weakly optimal, and the
condition of the theorem is best monotone.

\bigskip
\textbf{Proof of Theorem \ref{thm:def}:}\quad Suppose $\pi$ satisfies the
condition in Theorem \ref{thm:def}, but $\defi(G)\geq\beta+2$. (The
condition $\beta\equiv n\pmod2$ guarantees that $\defi(G)-\beta$ is always
even.) Define $G'\doteq K_{\beta+1}+G$, with degree sequence
$\pi'=(d_1+\beta+1,\dots,$ $d_n+\beta+1,((n-1)+\beta+1)^{\beta+1})$. Note
that the number of vertices of~$G'$ is odd.

Suppose $G'$ has a Hamilton cycle. Then, by taking alternating edges on
that cycle, there is a matching covering all vertices of $G'$ except one
vertex, and we can choose that missed vertex freely. So choose a matching
covering all but one of the $\beta+1$ new vertices. Removing the
other~$\beta$ new vertices as well, the remaining edges form a matching
covering all but at most~$\beta$ vertices from~$G$, a contradiction.

Hence~$G'$ cannot have a Hamilton cycle, and $\pi'$ cannot satisfy the
condition in Theorem \ref{thm:chvatal}. Thus there is some $i\ge\beta+1$
such that
\[d_i+\beta+1\le i<\tfrac12(n+\beta+1)\qquad\text{and}\qquad
d_{n+\beta+1-i}+\beta+1\le(n+\beta+1)-i-1.\]
Subtracting $\beta+1$ throughout this equation gives
$$d_i\le i-\beta-1<\tfrac12(n-\beta-1)\qquad\text{and}\qquad
d_{n+\beta+1-i}\le n-i-1.$$
Replacing $i$ by $j+1$ we get
$$d_{j+1}\le j-\beta<\tfrac12(n-\beta-1)\qquad\text{and}\qquad
d_{n+\beta-j}\le n-j-2.$$

Thus $\pi$ fails to satisfy the condition in Theorem \ref{thm:def}, a
contradiction.\eop

\bigskip
As an important special case, we give the best monotone condition for a
graph to have a 1-factor.

\begin{cor}[\cite{BC76,LV72}]\label{thm:1-factor}
  \; Let $G$ have degree sequence $\pi = (d_1 \leq \dots \leq d_n)$, with
  $n \geq 2$ and~$n$ even. If
  \begin{equation}\label{eqn0}
    d_{i+1} \leq i < \tfrac12 n\;\Longrightarrow\; d_{n-i} \geq n-i-1,
  \end{equation}
  then $G$ contains a 1-factor.
\end{cor}

We note in passing that \eqref{eqn0} is Chv\'atal's best monotone condition
for~$G$ to have a hamiltonian path~\cite{Chvatal72}.

\section{Best monotone condition for a 2-factor}\label{s3}

We now give a best monotone condition for the existence of a 2-factor. In
what follows we abuse the notation by setting $d_0=0$.

\begin{thm}\label{thm:2-factor}
  \quad Let $G$ have degree sequence $\pi = (d_1 \leq \dots \leq d_n)$,
  with $n \geq 3$. If{

    \qitemm{(i)} $n$ odd \ $\Longrightarrow$ \
    $d_{(n+1)/2}\ge\frac12(n+1)$;

    \qitemm{(ii)} $n$ even \ $\Longrightarrow$ \ $d_{(n-2)/2}\ge\frac12n$
    or $d_{(n+2)/2}\ge\frac12(n+2)$;

    \qitemm{(iii)} $d_i\le i$ and $d_{i+1}\le i+1$ \ $\Longrightarrow$ \
    $d_{n-i-1}\ge n-i-1$ or $d_{n-i}\ge n-i$, for $0\le i\le\frac12(n-2)$;

    \qitemm{(iv)} $d_{i-1}\le i$ and $d_{i+2}\le i+1$ \ $\Longrightarrow$ \
    $d_{n-i-3}\ge n-i-2$ or $d_{n-i}\ge n-i-1$, for $1\le
    i\le\frac12(n-5)$,

  }then $G$ contains a 2-factor.
\end{thm}

The condition in Theorem \ref{thm:2-factor} is easily seen to be monotone.
Furthermore, if $\pi$ fails to satisfy any of (i) through (iv), then $\pi$
is majorized by some $\pi'$ having a realization~$G'$ without a 2-factor.
In particular, note that

$\bullet$\quad if (i) fails, then $\pi$ is majorized by
$\pi'=\bigl(\frac12(n-1)\bigr)^{(n+1)/2}(n-1)^{(n-1)/2}$, having
realization $K_{(n-1)/2}+\overline{K_{(n+1)/2}}$;

$\bullet$\quad if (ii) fails, then $\pi$ is majorized by
$\pi'=\bigl(\frac12(n-2)\bigr)^{(n-2)/2}\bigl(\frac12n\bigr)^2
(n-1)^{(n-2)/2}$, having realization
$K_{(n-2)/2}+(\overline{K_{(n-2)/2}}\cup K_2)$;

$\bullet$\quad if (iii) fails for some $i$, then $\pi$ is majorized by
$\pi'=i^i(i+1)^1(n-i-2)^{n-2i-2}$ $(n-i-1)^1(n-1)^i$, having realization
$K_i+(\overline{K_{i+1}}\cup K_{n-2i-1})$ together with an edge joining
$\overline{K_{i+1}}$ and $K_{n-2i-1}$;

$\bullet$\quad if (iv) fails for some $i$, then $\pi$ is majorized by
$\pi'= i^{i-1}(i+1)^3(n-i-3)^{n-2i-5}$ $(n-i-2)^3(n-1)^i$, having
realization $K_i + (\overline{K_{i+2}} \cup K_{n-2i-2})$ together with
three independent edges joining $\overline{K_{i+2}}$ and $K_{n-2i-2}$.

It is immediate that none of the above realizations contain a 2-factor.
Hence, Theorem \ref{thm:2-factor} is weakly optimal, and the condition of
the theorem is best monotone.

\bigskip
\textbf{Proof of Theorem \ref{thm:2-factor}:}\quad Suppose $\pi$ satisfies
(i) through (iv), but $G$ has no \mbox{2-factor}. We may assume the
addition of any missing edge to $G$ creates a 2-factor. Let
$v_1, \dots, v_n$ be the vertices of $G$, with respective degrees
$d_1\leq \dots \leq d_n$, and assume $v_j, v_k$ are a nonadjacent pair with
$j+k$ as large as possible, and $d_j \leq d_k$. Then~$v_j$ must be adjacent
to $v_{k+1}, v_{k+2}, \dots, v_n$ and so
\begin{equation}\label{eqn6}
  d_j \geq n-k.
\end{equation}
Similarly, $v_k$ must be adjacent to
$v_{j+1}, \dots, v_{k-1}, v_{k+1}, \dots,v_n$, and so
\begin{equation}\label{eqn7}
  d_k \geq n-j-1.
\end{equation}
Since $G+(v_j,v_k)$ has a 2-factor, $G$ has a spanning subgraph consisting
of a path~$P$ joining $v_j$ and $v_k$, and $t \geq 0$ cycles
$C_1,\dots,C_t$, all vertex disjoint.

We may also assume $v_j,v_k$ and $P$ are chosen such that if $v,w$ are any
nonadjacent vertices with $d_G(v)=d_j$ and $d_G(w)=d_k$, and if $P'$ is any
$(v,w)$-path such that $G-V(P')$ has a \mbox{2-factor}, then $|P'|\le|P|$.
Otherwise, re-index the set of vertices of degree $d_j$ (resp., $d_k$) so
that $v$ (resp., $w$) is given the highest index in the set.

Since $G$ has no \mbox{2-factor}, we cannot have independent edges between
$\{v_j,v_k\}$ and two consecutive vertices on any of the $C_{\mu}$,
$0\le\mu\le t$. Similarly, we cannot have $d_P(v_j)+d_P(v_k) \ge |V(P)|$,
since otherwise $\langle V(P) \rangle$ is hamiltonian and $G$ contains a
\mbox{2-factor}. This means
\begin{equation}\label{eq7b}
  \begin{array}{@{}c@{}}
    d_{C_{\mu}}(v_j) + d_{C_{\mu}}(v_k) \leq |V(C_{\mu})|\quad\text{for
      $0 \leq \mu \leq t$},\\[4pt]
    \text{and}\quad d_P(v_j) + d_P(v_k) \leq |V(P)|-1.
  \end{array}
\end{equation}
It follows immediately that
\begin{equation}\label{eq7c}
  d_j + d_k \leq n-1.
\end{equation}
We distinguish two cases for $d_j+d_k$.

\bigskip
\textsc{Case 1:}\quad $d_j + d_k \leq n-2$.

\smallskip
Using \eqref{eqn7}, we obtain
$$d_j \leq (n-2) - d_k \leq (n-2) - (n-j-1) = j-1.$$
Take $i,m$ so that $i=d_j=j-m$, where $m\ge1$. By Case~1 we have
$i\le\frac12(n-2)$. Since also $d_i=d_{j-m}\le d_j=i$ and
$d_{i+1}=d_{j-m+1}\le d_j=i$, condition~(iii) implies
$d_{n-(j-m)-1} \geq n-(j-m)-1$ or $d_{n-(j-m)} \geq n-(j-m)$. In either
case,
\begin{equation}\label{eqn9}
  d_{n-(j-m)} \geq n-(j-m)-1.
\end{equation}
Adding $d_j = j-m$ to \eqref{eqn9}, we obtain
\begin{equation}\label{eqn10}
  d_j + d_{n-j+m} \geq n-1.
\end{equation}
But $d_j+d_k \leq n-2$ and \eqref{eqn10} together give $n-j+m>k$, hence
$j+k<n+m$. On the other hand, \eqref{eqn6} gives $j-m=d_j \geq n-k$, hence
$j+k \geq n+m$, a contradiction.\openeop

\bigskip
\textsc{Case 2:}\quad $d_j + d_k = n-1$.

\smallskip
In this case we have equality in~\eqref{eq7c}, hence all the inequalities
in~\eqref{eq7b} become equalities. In particular, this implies that every
cycle $C_{\mu}$, $1 \leq \mu \leq t$, satisfies one of the following
conditions:{

  \qitemm{(a)} Every vertex in $C_{\mu}$ is adjacent to $v_j$ (resp.,
  $v_k$), and none are adjacent to $v_k$ (resp., $v_j$), or

  \qitemm{(b)}$|V(C_{\mu})|$ is even, and $v_j, v_k$ are both adjacent to
  the same alternate vertices on~$C_{\mu}$.

}We call a cycle of type~(a) a \textit{$j$-cycle} (resp.,
\textit{$k$-cycle}), and a cycle of type~(b) a \textit{$(j,k)$-cycle}. Set
$A = \bigcup_{\textrm{$j$-cycles $C$}}V(C)$,
$B =\bigcup_{\textrm{$k$-cycles $C$}} V(C)$, and
$D=\bigcup_{\textrm{$(j,k)$-cycles $C$}}V(C)$, and let $a \doteq |A|$,
$b\doteq |B|$, and $c\doteq\frac12|D|$.

Vertices in $V(G)-\{v_j,v_k\}$ which are adjacent to both (resp., neither)
of $v_j,v_k$ will be called \textit{large} (resp., \textit{small})
vertices. In particular, the vertices of each $(j,k)$-cycle are alternately
large and small, and hence there are $c$ small and $c$ large vertices among
the $(j,k)$-cycles.

By the definitions of $a,b,c$, noting that a cycle has at least 3 vertices,
we have the following.

\begin{obs}\label{obs1}
  \quad We have $a=0$ or $a\ge3$, $b=0$ or $b\ge3$, and $c=0$ or $c\ge2$.
\end{obs}

By the choice of~$v_j,v_k$ and~$P$, we also have the following
observations.

\begin{obs}\label{obs2}
\mbox{}

 \vspace*{-\parskip}
 \qitemm{(a)}If $(u,v_k) \notin E(G)$, then $d_G(u)\le d_j$; if
 $(u,v_j)\notin E(G)$, then $d_G(u)\le d_k$.

 \qitemm{(b)}A vertex in~$A$ has degree at most $d_j-1$.

 \qitemm{(c)}A vertex in~$B$ has degree at most $d_k-1$.

 \qitemm{(d)}A small vertex in~$D$ has degree at most $d_j-1$.
\end{obs}

\textbf{Proof:}\quad Part~(a) follows directly from the choice of $v_j,v_k$
as nonadjacent with $d_G(v_j)+d_G(v_k)=d_j+d_k$ maximal.

For~(b), consider any $a \in A$, with say $a \doteq v_{\ell}$. Since
$(v_{\ell},v_k) \notin E(G)$, we have $\ell<j$ by the maximality of $j+k$,
and so $d_G(a) \le d_j$. If $d_G(a)=d_j$, then since each vertex in $A$ is
adjacent to $v_j$, we can combine the path $P$ and the $j$-cycle $C_{\mu}$
containing $a$ (leaving the other cycles $C_{\mu}$ alone) into a path $P'$
joining $a$ and $v_k$ such that $G-V(P')$ has a \mbox{2-factor} and
$|P'|>|P|$, contradicting the choice of $P$. Thus $d_G(a) \le d_j-1$,
proving~(b).

Parts~(c) and~(d) follow by a similar arguments.\openeop

\bigskip
Let $p\doteq|V(P)|$, and let us re-index $P$ as
$v_j=w_1,w_2,\dots,w_p=v_k$. By the case assumption,
$d_P(w_1)+d_P(w_p)=p-1$.

Assume first that $p=3$. Then $d_j = a+c+1$ and $d_k=b+c+1$, so that
$b \ge a$. Moreover, $n=a+b+2c+3$ and there are $c+1$ large vertices and
$c$ small vertices.

If $b \ge 3$, the large vertex $w_2$ is not adjacent to a vertex in $A$ or
to a small vertex in~$D$, or else $G$ contains a 2-factor. Thus $w_2$ has
degree at most $n-1-(a+c)$, and by Observations 2(b,c,d), $\pi$ is
majorized by
\[\pi_1 = (a+c)^{a+c} (a+c+1)^1 (b+c)^b (b+c+1)^1 (n-1-(a+c))^1 (n-1)^c.\]
Setting $i=a+c$, so that $0 \le i = a+c = (n-3)-(b+c) \le \frac12(n-3)$,
$\pi_1$ becomes
\[\pi_1 = i^i (i+1)^1 (n-i-3)^b (n-i-2)^1 (n-i-1)^1 (n-1)^c.\]
Since $\pi_1$ majorizes $\pi$, we have $d_i \le i$, $d_{i+1} \le i+1$,
$d_{n-i-1} = d_{n-(a+c+1)} \le n-i-2$, and $d_{n-i}=d_{n-(a+c)} \le n-i-1$,
and $\pi$ violates condition~(iii). Hence $b=0$ by Observation 1, and a
fortiori $a=0$.

But if $a=b=0$, then $c=\frac12(n-3)$, $n$ is odd, and by Observation 2(d),
$\pi$ is majorized by
\[\pi_2 = \bigl(\tfrac12(n-3)\bigr)^{(n-3)/2} \bigl(\tfrac12(n-1)\bigr)^2
(n-1)^{(n-1)/2}.\]
Since $\pi_2$ majorizes $\pi$, we have $d_{(n+1)/2} \leq \frac12(n-1)$, and
$\pi$ violates condition~(i).

Hence we assume $p \ge 4$.

We make several further observations regarding the possible adjacencies of
$v_j,v_k$ into the path $P$.

\begin{obs}\label{obs3}
  \quad For all $m$, $1 \leq m \leq p-1$, we have $(w_1, w_{m+1}) \in E(G)$
  if and only if $(w_p, w_m) \notin E(G)$.
\end{obs}

\textbf{Proof:}\quad If $(w_1,w_{m+1}) \in E(G)$ then,
$(w_p,w_m) \notin E(G)$, since otherwise $\langle V(P) \rangle$ is
hamiltonian and $G$ has a 2-factor. The converse follows since
$d_P(w_1)+d_P(w_p) = p-1$.\openeop

\begin{obs}\label{obs4}
  \quad If $(w_1,w_m),(w_1,w_{m+1}) \in E(G)$ for some $m$,
  $3 \leq m \leq p-3$, then we have $(w_1,w_{m+2}) \in E(G)$.
\end{obs}

\textbf{Proof:}\quad If $(w_1,w_{m+2}) \notin E(G)$, then
$(w_p,w_{m+1}) \in E(G)$ by Observation~\ref{obs3}. But since
$(w_1,w_m)\in E(G)$, this means that $\langle V(P) \rangle$ would have a
2-factor consisting of the cycles $(w_1, w_2, \dots,w_m, w_1)$ and
$(w_p,w_{m+1},w_{m+2}, \dots, w_p)$, and thus $G$ would have a
\mbox{2-factor}, a contradiction.\openeop

\bigskip
Observation~\ref{obs4} implies that if $w_1$ is adjacent to consecutive
vertices $w_m,w_{m+1} \in V(P)$ for some $m \geq 3$, then $w_1$ is adjacent
to all of the vertices $w_m, w_{m+1}, \dots, w_{p-1}$.

\begin{obs}\label{obs5}
  \quad If $(w_1,w_m),(w_1,w_{m-1}) \notin E(G)$ for some
  $5 \leq m \leq p-1$, then we have $(w_1,w_{m-2}) \notin E(G)$.
\end{obs}

\textbf{Proof:}\quad If $(w_1,w_m) \notin E(G)$, then $(w_p,w_{m-1}) \in
E(G)$ by Observation~\ref{obs3}. So if also $(w_1,w_{m-2})\in E(G)$, then
$\langle V(P) \rangle$ would have a 2-factor as in the proof of
Observation~\ref{obs4}, leading to the same contradiction.\openeop

\bigskip
Observation~\ref{obs5} implies that if $w_1$ is not adjacent to two
consecutive vertices $w_{m-1}, w_m$ on $P$ for some $m \leq p-1$, then
$w_1$ is not adjacent to any of $w_3,\ldots,w_{m-1},w_m$.

\medskip
By Observation~\ref{obs3}, the adjacencies of $w_1$ into $P$ completely
determine the adjacencies of $w_p$ into $P$. But combining Observations
\ref{obs4} and \ref{obs5}, we see that the adjacencies of~$w_1$ and $w_p$
into $P$ must appear as shown in Figure 1, for some $\ell, r \geq 0$. In
summary, $w_1$ will be adjacent to $r \geq 0$ consecutive vertices
$w_{p-r}, \dots, w_{p-1}$ (where $w_{\alpha}, \dots, w_{\beta}$ is taken to
be empty if $\alpha>\beta$), $w_p$ will be adjacent to $\ell \geq 0$
consecutive vertices $w_2,\dots,w_{\ell+1}$, and $w_1,w_p$ are each
adjacent to the vertices $w_{\ell+3}, w_{\ell+5}, \dots, w_{p-r-4},
w_{p-r-2}$. Note that $\ell=p-2$ implies $r=0$, and $r=p-2$ implies
$\ell=0$.

\begin{figure}[ht]
\centering
\includegraphics[width=0.6\textwidth]{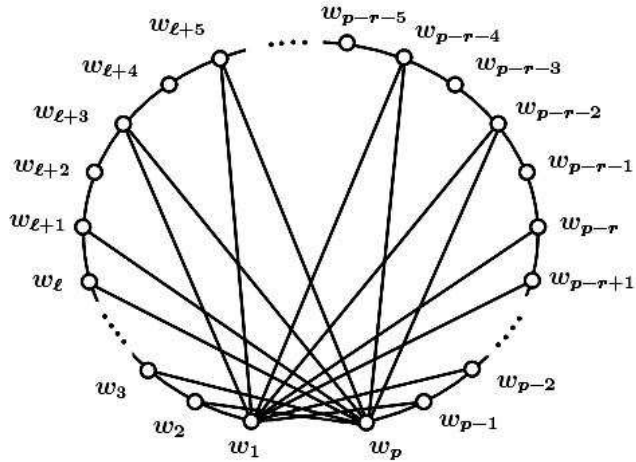}
\caption{The adjacencies of $w_1, w_p$ on $P$.}
\end{figure}

Counting neighbors of~$w_1$ and~$w_p$ we get their degrees as follows.

\begin{obs}\label{obs6}
  \begin{align*}
    d_j=d_G(w_1)&=\left\{\begin{array}{ll}
        a+c+1,&\text{if $\ell=p-2$, $r=0$},\\[4pt]
        a+c+p-2,&\text{if $r=p-2$, $\ell=0$},\\[4pt]
        a+c+r+\tfrac12(p-r-\ell-1);&\text{otherwise};\end{array}\right.\\[2mm]
    d_k=d_G(w_p)&=\left\{\begin{array}{ll}
        b+c+p-2,&\text{if $\ell=p-2$, $r=0$},\\[4pt]
        b+c+1,&\text{if $r=p-2$, $\ell=0$},\\[4pt]
        b+c+\ell+\tfrac12(p-r-\ell-1);&\text{otherwise}.\end{array}\right.
  \end{align*}
\end{obs}

We next prove some observations to limit the possibilities for $(a,b)$ and
$(\ell,r)$.

\begin{obs}\label{obs7}
  \quad If $(w_1,w_{p-1}) \in E(G)$ (resp., $(w_2,w_p) \in E(G))$, then we
  have $b=0$ (resp., $a=0$).
\end{obs}

\textbf{Proof:}\quad If $b \neq 0$, there exists a $k$-cycle
$C\doteq(x_1, x_2, \dots, x_s, x_1)$. But if also $(w_1,w_{p-1}) \in E(G)$,
then $(w_1,w_2, \dots, w_{p-1}, w_1)$ and $(w_p, x_1, \dots, x_s, w_p)$
would be a 2-factor in $\langle V(C) \cup V(P) \rangle$, implying a
2-factor in $G$. The proof that $(w_2,w_p) \in E(G)$ implies $a=0$ is
symmetric.\openeop

\bigskip
From Observation \ref{obs6}, we have
\begin{equation}\label{eqstar}
  0 \le d_k-d_j = b-a + \left\{\begin{array}{ll}
      p-3,&\text{if $\ell=p-2$, $r=0$},\\[4pt]
      3-p,&\text{if $r=p-2$, $\ell=0$},\\[4pt]
      \ell-r,&\text{otherwise}.\end{array}\right.\\[2mm]
\end{equation}

From this, we obtain

\begin{obs}\label{obs8}
  \quad $\ell \geq r$.
\end{obs}

\textbf{Proof:}\quad Suppose first $r \ne p-2$. If $r > \ell \ge 0$, then
$b > a \ge 0$ since $b + \ell \ge a+r$ by \eqref{eqstar}. But $r>0$ implies
$(w_1,w_{p-1}) \in E(G)$, and thus $b=0$ by Observation \ref{obs7}, a
contradiction.

Suppose then $r = p-2 \ge 2$. Then $b > a \ge 0$, since $b \ge a+p-3$ by
\eqref{eqstar}. Since $r>0$, we have the same contradiction as in the
previous paragraph.\openeop

\begin{obs}\label{obs9}
  \quad If $r \geq 1$, then $\ell \leq 1$.
\end{obs}

\textbf{Proof:}\quad Else we have
$(w_1, w_{p-1}), (w_p, w_2), (w_p,w_3) \in E(G)$, and $(w_1, w_2, w_p, w_3,
\dots,\linebreak w_{p-1}, w_1)$ would be a hamiltonian cycle in $\langle
V(P)\rangle$. Thus $G$ would have a 2-factor, a contradiction.\openeop

\bigskip
Observations~\ref{obs8} and~\ref{obs9} together limit the possibilities for
$(\ell, r)$ to $(1,1)$ and $(\ell,0)$ with $0\le\ell\le p-2$. We also
cannot have $(\ell,r)=(p-3,0)$, since~$w_p$ is always adjacent to
$w_{p-1}$, and so we would have $\ell=p-2$ in that case. And we cannot have
$(\ell,r)=(p-4,0)$, since then $p-r-\ell-1$ is odd, violating
Observation~\ref{obs6}. To complete the proof of
Theorem~\ref{thm:2-factor}, we will deal with the remaining possibilities
in a number of cases, and show that all of them lead to a contradiction of
one or more of conditions (i) through (iv).

Before doing so, let us define the spanning subgraph $H$ of $G$ by letting
$E(H)$ consist of the edges in the cycles $C_{\mu}$, $0 \le \mu \le t$, or
in the path $P$, together with the edges incident to $w_1$ or $w_p$. Note
that the edges incident to $w_1$ or $w_p$ completely determine the large or
small vertices in $G$. In the proofs of the cases below, any adjacency
beyond those indicated would create an edge $e$ such that $H+e$, and a
fortiori $G$, contains a \mbox{2-factor}.

\bigskip
\textsc{Case 2.1:}\quad $(\ell, r) = (1,1)$.

\smallskip
Since $(w_1,w_{p-1}), (w_2,w_p) \in E(G)$, we have $a=b=0$, by
Observation~\ref{obs7}. Using Observation~\ref{obs6} this means that
$d_j=d_k=\frac12(n-1)$, and hence~$n$ is odd. Additionally, there are
$c+\frac12(p-3)=\frac12(n-3)$ small vertices. Each of these small vertices
has degree at most~$d_j$ by Observation~\ref{obs2}\,(a), and so $\pi$ is
majorized by
\[\pi_3=\bigl(\tfrac12(n-1)\bigr)^{(n+1)/2}(n-1)^{(n-1)/2}.\]
But~$\pi_3$ (a fortiori $\pi$) violates condition (i).\openeop

\bigskip
\textsc{Case 2.2:}\quad $(\ell, r)= (0,0)$.

\smallskip
By Observation~\ref{obs6}, $d_j =a+ c+ \frac12(p-1)$ and
$d_k =b+ c+ \frac12(p-1)$, so that $b \geq a$. Also, there are
$c+\frac12(p-3)$ large and $c+\frac12(p-5)$ small vertices.

$\bullet$\quad By Observation~\ref{obs2}\,(b,c), each vertex in~$A$
(resp.,~$B$) has degree at most $d_j-1=a+ c+ \frac12(p-3)$ (resp.,
$d_k-1=b+ c+ \frac12(p-3)$).

$\bullet$\quad Each small vertex is adjacent to at most the large vertices
(otherwise~$G$ contains a 2-factor), and so each small vertex has degree at
most $c+\frac12(p-3)$.

$\bullet$\quad The vertex $w_2$ (resp., $w_{p-1}$) is adjacent to at most
the large vertices and $w_1$ (resp., $w_p$) (otherwise~$G$ contains a
2-factor), and so $w_2$, $w_{p-1}$ each have degree at most
$c+\frac12(p-1)$.

Thus~$\pi$ is majorized by
\begin{align*}
 \pi_4={}&\bigl(c+\tfrac12(p-3)\bigr)^{c+(p-5)/2}
 \bigl(c+\tfrac12(p-1)\bigr)^2\bigl(a+c+\tfrac12(p-3)\bigr)^a\\
 &\qquad\bigl(a+c+\tfrac12(p-1)\bigr)^1
 \bigl(b+c+\tfrac12(p-3)\bigr)^b\bigl(b+c+\tfrac12(p-1)\bigr)^1
 (n-1)^{c+(p-3)/2}.
\end{align*}

Setting $i=a+c+\frac12(p-1)$, so that
$2 \le i = \frac12(n-(b-a)-1) \le \frac12(n-1)$, the sequence $\pi_4$
becomes
\[\pi_4=(i-a-1)^{i-a-2} (i-a)^2 (i-1)^a i^1 (n-i-2)^{n-2i+a-1}
(n-i-1)^1 (n-1)^{i-a-1}.\]

If $2 \le i \le \frac12(n-2)$, then since $\pi_4$ majorizes $\pi$, we have
$d_i \le i$, $d_{i+1} \le i$, $d_{n-i-1} \le n-i-2$, and
$d_{n-i} \le n-i-2$, and $\pi$ violates condition~(iii).

If $i=\frac12(n-1)$, then $n$ is odd, and $\pi_4$ reduces to
\[\pi_4' = \bigl(\tfrac12(n-3)-a\bigr)^{(n-5)/2-a}
\bigl(\tfrac12(n-1)-a\bigr)^2 \bigl(\tfrac12(n-3)\bigr)^{2a}
\bigl(\tfrac12(n-1)\bigr)^2 (n-1)^{(n-3)/2-a}.\]
Since $\pi_4'$ majorizes $\pi$, we have $d_{(n+1)/2} \le \frac12(n-1)$, and
$\pi$ violates condition~(i).\openeop

\bigskip
\textsc{Case 2.3:}\quad $(\ell, r) = (1,0)$

\smallskip
By Observation~\ref{obs7}, $a=0$, and thus by Observation \ref{obs6},
$d_j=c+\frac12(p-2)$ and $d_k=b+c+\frac12p$. Also, there are
$c+\frac12(p-2)$ large and $c+\frac12(p-4)$ small vertices. If $p=4$ then
$\ell=2$, a contradiction, and hence $p\ge6$.

$\bullet$\quad By Observation~\ref{obs2}\,(c), each vertex in~$B$ has
degree at most $d_k-1=b+c+\frac12(p-2)$.

$\bullet$\quad Each small vertex is adjacent to at most the large vertices, and so each small vertex has degree at most $c+\frac12(p-2)$.

$\bullet$\quad The vertex $w_{p-1}$ is adjacent to at most~$w_p$ and the
large vertices, and so $w_{p-1}$ has degree at most $c+\frac12p$.

Thus~$\pi$ is majorized by
\[\pi_5=\bigl(c+\tfrac12(p-2)\bigr)^{c+(p-2)/2}
 \bigl(c+\tfrac12p\bigr)^1\bigl(b+c+\tfrac12(p-2)\bigr)^b
\bigl(b+c+\tfrac12p\bigr)^1(n-1)^{c+(p-2)/2}.\]
Setting $i=c+\frac12(p-2)$, so that
$2 \le i = \frac12(n-b-2) \le \frac12(n-2)$, $\pi_5$ becomes
\[\pi_5 = i^i (i+1)^1(n-i-2)^{n-2i-2}(n-i-1)^1(n-1)^i.\]
If $2 \le i \le \frac12(n-3)$, then since $\pi_5$ majorizes $\pi$, we have
$d_i \le i$, $d_{i+1} \le i+1$, $d_{n-i-1} \le n-i-2$, and
$d_{n-i} \le n-i-1$, and $\pi$ violates condition~(iii).

If $i = \frac12(n-2)$, then $n$ is even, and $\pi_5$ reduces to
\[\pi_5'= \bigl(\tfrac12n-1\bigr)^{n/2-1}
\bigl(\tfrac12n\bigr)^2(n-1)^{n/2-1}.\]
Since $\pi_5'$ majorizes $\pi$, we have $d_{n/2-1} \le \tfrac12n-1$ and
$d_{n/2+1}\le\tfrac12n$, and $\pi$ violates condition~(ii).\openeop

\bigskip
\textsc{Case 2.4:}\quad $(\ell, r) = (\ell, 0)$, \textit{where}
$2 \leq \ell\leq p-5$

\smallskip
We have $a=0$ by Observation \ref{obs7}, and $p-\ell\ge 5$ by Case 2.4. By
Observation~\ref{obs6}, $d_j=c+\frac12(p-\ell-1)$ and
$d_k=b+c+\ell+\frac12(p-\ell-1)$. Moreover, there are $c+\frac12(p-\ell-1)$
large vertices including~$w_2$, and $c+\frac12(p-\ell-3)$ small vertices.

$\bullet$\quad By Observation~\ref{obs2}\,(c), each vertex in~$B$ has
degree at most $d_k-1=b+c+\ell+\frac12(p-\ell-3)$.

$\bullet$\quad Each small vertex other than $w_{\ell+2}$ is adjacent to at
most the large vertices except $w_2$, and so each small vertex other than
$w_{\ell+2}$ has degree at most $c+\frac12(p-\ell-3)$.

$\bullet$\quad The vertex $w_{\ell+2}$ is not adjacent to $w_p$, and so by
Observation \ref{obs2}\,(a), $w_{\ell+2}$ has degree at most
$d_j=c+\frac12(p-\ell-1)$.

$\bullet$\quad The vertex $w_{p-1}$ is adjacent to at most $w_p$ and the
large vertices except $w_2$, and so $w_{p-1}$ has degree at most
$c+\frac12(p-\ell-1)$.

$\bullet$\quad Each $w_m$, $3 \le m \le \ell$, is adjacent to at most $w_p$, the large vertices, the vertices in~$B$, and $\{w_3,\dots,w_{\ell+1}\}-\{w_m\}$.  Hence each such $w_m$ has degree at most $b+c+\ell+\frac12(p-\ell-3)$.

$\bullet$\quad The vertex $w_2$ is adjacent to at most $w_1$, $w_p$, the
other large vertices, the vertices in $B$, and $\{w_3,\dots,w_{\ell+1}\}$.
Hence $w_2$ has degree at most $b+c+\ell+\frac12(p-\ell-1)$.

$\bullet$\quad The vertex $w_{\ell+1}$ is not adjacent to $w_1$, and so by
Observation \ref{obs2}\,(a), vertex $w_{\ell+1}$ has degree at most
$d_k=b+c+\ell+\frac12(p-\ell-1)$.

Thus~$\pi$ is majorized by
\begin{align*}
 &\qquad\pi_6=\bigl(c+\tfrac12(p-\ell-3)\bigr)^{c+(p-\ell-5)/2}
 \bigl(c+\tfrac12(p-\ell-1)\bigr)^3\\
 &\bigl(b+c+\ell+\tfrac12(p-\ell-3)\bigr)^{b+\ell-2}
 \bigl(b+c+\ell+\tfrac12(p-\ell-1)\bigr)^3(n-1)^{c+(p-\ell-3)/2}.
\end{align*}
Setting $i=c-1+\frac12(p-\ell-1)$, so that
$1 \le i = \frac12(n-b-\ell-3) \le \frac12(n-5)$, $\pi_6$ becomes
$$\pi_6=i^{i-1}(i+1)^3(i+b+\ell)^{b+\ell-2}(i+b+\ell+1)^3(n-1)^i.$$
Since $\pi_6$ majorizes $\pi$, we have $d_{i-1} \le i$, $d_{i+2} \le i+1$,
$d_{n-i-3} \le i+b+\ell = n-i-3$, and $d_{n-i} \le i+b+\ell+1 = n-i-2$, and
thus $\pi$ violates condition~(iv).\openeop

\bigskip
\textsc{Case 2.5:}\quad $(\ell, r) = (p-2,0)$

\smallskip
We have $a=0$, by Observation \ref{obs7}. By Observation \ref{obs6}, we
then have $d_j=c+1$ and $d_k = b+c+p-2$. If $d_1\le1$, then condition~(iii)
with $i=0$ implies $d_{n-1}\ge n-1$, which means there are at least~2
vertices adjacent to all other vertices, a contradiction. Hence $c+1=d_j\ge
d_1\ge2$, and so
$c \ge 2$ by Observation 1. Finally, there are $c+1$ large vertices
including $w_2$, and $c$ small vertices.

$\bullet$\quad By Observation~\ref{obs2}\,(a), the vertices in~$B$ have
degree at most $d_k=b+c+p-2$.

$\bullet$\quad By Observation~\ref{obs2}\,(d), the small vertices in~$D$
have degree at most $d_j-1=c$.

$\bullet$\quad The vertex $w_2$ is not adjacent to the small vertices in
$D$, and so $w_2$ has degree at most $n-1-c=b+c+p-1$.

$\bullet$\quad The vertices $w_3,\ldots,w_{p-1}$ have degree at most
$d_k=b+c+p-2$ by Observation~\ref{obs2}\,(a), since none of them are
adjacent to $w_1=v_j$.

Thus~$\pi$ is majorized by
\[\pi_7=c^c(c+1)^1(b+c+p-2)^{b+p-2}(b+c+p-1)^1(n-1)^c.\]
Setting $i=c$, so that $2 \le c = i = \frac12(n-b-p) \le \frac12(n-4)$,
$\pi_7$ becomes
\[\pi_7=i^i(i+1)^1(n-i-2)^{n-2i-2}(n-i-1)^1(n-1)^i.\]
Since $\pi_7$ majorizes $\pi$, we have $d_i \le i$, $d_{i+1} \le i+1$,
$d_{n-i-1} \le n-i-2$, and $d_{n-i} \le n-i-1$, and $\pi$ violates
condition~(iii).\openeop

\bigskip
The proof of Theorem \ref{thm:2-factor} is complete.\eop

\section{Sufficient condition for the existence of a\\
  {\boldmath$k$-factor, $k \geq 2$}}\label{s4}

The increase in complexity of Theorem \ref{thm:2-factor} ($k=2$) compared
to Corollay \ref{thm:1-factor} ($k=1$) suggests that the best monotone
condition for $\pi$ to be forcibly $k$-factor graphical may become unwieldy
as $k$ increases. Indeed, we make the following conjecture.

\begin{conj}\label{conj:superpoly}
  \quad The best monotone condition for a degree sequence $\pi$ of length
  $n$ to be forcibly $k$-factor graphical requires checking at least $f(k)$
  nonredundant conditions (where each condition may require $O(n)$ checks),
  where $f(k)$ grows superpolynomially in $k$.
\end{conj}

Kriesell \cite{Kp} has verified such rapidly increasing complexity for the
best monotone condition for~$\pi$ to be forcibly $k$-edge-connected.
Indeed, Kriesell has shown such a condition entails checking at least
$p(k)$ nonredundant conditions, where $p(k)$ denotes the number of
partitions of $k$. It is well-known \cite{HR18} that
$p(k) \sim \dfrac{e^{\pi \sqrt{2k/3}}}{4 \sqrt{3} k} $.

The above conjecture suggests the desirability of obtaining a monotone
condition for~$\pi$ to be forcibly $k$-factor graphical which does not
require checking a superpolynomial number of conditions. Our goal in this
section is to prove such a condition for $k \geq 2$. Since our condition
will require Tutte's Factor Theorem \cite{B50, T52}, we begin with some
needed background.

Belck \cite{B50} and Tutte \cite{T52} characterized graphs $G$ that do not
contain a $k$-factor. For disjoint subsets $A,B$ of $V(G)$, let
$C=V(G)-A-B$. We call a component $H$ of $\langle C\rangle$ \emph{odd} if
$k|H| + e(H,B)$ is odd. The number of odd components of $\langle C\rangle$
is denoted by $odd_k(A,B)$. Define
$$\Theta_k(A,B) \doteq   k|A| + \sum_{u \in B} d_{G-A}(u) - k|B| - 
odd_k(A,B).$$

\begin{thm}\label{thm:5.2}
  \quad Let $G$ be a graph on $n$ vertices and $k \geq 1$.

  \qitemm{(a)} \emph{\cite{T52}.} For any disjoint $A,B \subseteq V(G)$,
  $\Theta_k(A,B) \equiv kn\pmod2$;

  \qitemm{(b)} \emph{\cite{B50,T52}.} The graph $G$ does not contain a
  $k$-factor if and only if $\Theta_k(A,B)<0$, for some disjoint $A,B
  \subseteq V(G)$.
\end{thm}

We call any disjoint pair $A,B \subseteq V(G)$ for which $\Theta_k(A,B)<0$
a \textit{$k$-Tutte-pair} for~$G$. Note that if $kn$ is even, then $A,B$ is
a $k$-Tutte-pair for~$G$ if and only if
\[k|A| + \sum_{u \in B} d_{G-A}(u)\le k|B| + odd_k(A,B) -2.\]
Moreover, for all $u\in B$ we have $d_G(u)\le d_{G-A}(u) +|A|$, so
$\sum\limits_{u\in B}d_G(u)\le\sum\limits_{u\in B}d_{G-A}(u)+|A||B|$. Thus
for each $k$-Tutte-pair $A,B$ we have
\begin{equation}\label{eqn12}
\sum_{u\in B}d_G(u)\le k|B|+|A||B|-k|A|+odd_k(A,B)-2.
\end{equation}

Our main result in this section is the following condition for a graphical
degree sequence $\pi$ to be forcibly $k$-factor graphical. The condition
will guarantee that no $k$-Tutte-pair can exist, and is readily seen to be
monotone. We again set $d_0=0$.

\begin{thm}\label{thm:k-factor}
  \quad Let $\pi = (d_1 \leq \dots \leq d_n)$ be a graphical degree
  sequence, and let $k \geq 2$ be an integer such that $kn$ is even.
  Suppose{

    \qitemm{(i)}$d_1\ge k$;

    \qitemm{(ii)}for all $a,b,q$ with $0\le a<\frac12 n$, $0\le b\le n-a$
    and $\max\{0,a(k-b)+2\}\le q\le n-a-b$ so that
    $\sum\limits_{i=1}^bd_i\le kb+ab-ka+q-2$, the following holds: Setting
    $r=a+k+q-2$ and $s=n-\max\{0,b-k+1\}-\max\{0,q-1\}-1$, we have

    \leavevmode\hbox to10mm{\hss}\hbox to10mm{$(*)$\hss}$r\le s$ and
    $d_b\le r$, or $r>s$ and $d_{n-a-b}\le s$ \ $\Longrightarrow$ \
    $d_{n-a}\ge\max\{r,s\}+1$.

  }Then $\pi$ is forcibly $k$-factor graphical.
\end{thm}

\textbf{Proof}:\quad Let $n$ and $k\ge2$ be integers with $kn$ even.
Suppose $\pi$ satisfies~(i) and~(ii) in the theorem, but has a
realization~$G$ with no $k$-factor. This means that~$G$ has at least one
$k$-Tutte-pair.

Following \cite{J4}, a $k$-Tutte-pair $A,B$ is \textit{minimal} if either
$B = \varnothing$, or $\Theta_k(A,B') \geq 0$ for all proper subsets $B'
\subset B$. We then have

\begin{lem}[\cite{J4}]\label{lem:J4lem}
  \; Let $k \geq 2$, and let $A,B$ be a minimal $k$-Tutte-pair for a
  graph~$G$ with no $k$-factor. If $B \neq \varnothing$, then
  $\Delta(\langle B \rangle) \leq k-2$.
\end{lem}

Next let $A,B$ be a \mbox{$k$-Tutte-pair} for $G$ with $A$ as large as
possible, and $A,B$ minimal. Also, set $C=V(G)-A-B$. We establish some
further observations.

\begin{lem}\label{lem:5.4}
\mbox{}

 \vspace*{-\parskip}
  \qitemm{(a)}$|A| < \frac12n$.

 \qitemm{(b)}For all $v \in C$, $e(v,B) \leq \min\{k-1,|B|\}$.

 \qitemm{(c)}For all $u \in B$, $d_G(u) \leq|A|+ k + odd_k(A,B)-2$.
\end{lem}

\textbf{Proof:}\quad Suppose $|A| \geq \frac12 n$, so that $|A|\ge|B|+|C|$.
Then we have
\begin{align*}
  \Theta_k(A,B)&{}=k|A|+\sum_{u \in B} d_{G-A}(u)-k|B|-odd_k(A,B)\geq
  k(|A|-|B|)-odd_k(A,B)\\
  &{}\ge k|C|-odd_k(A,B)>|C|-odd_k(A,B)\ge0,
\end{align*}
which contradicts that $A,B$ is a $k$-Tutte-pair.

For (b), clearly $e(v,B) \leq |B|$. If $e(v,B) \geq k$ for some $v \in C$,
move $v$ to $A$, and consider the change in each term in $\Theta_k(A,B)$:
$$\underbrace{k|A|\rule[-5mm]{0pt}{0pt}}_{\textrm{increases by $k$}} + 
\underbrace{{\textstyle\sum\limits_{u \in B}}
  d_{G-A}(u)\rule[-5mm]{0pt}{0pt}}_{\textrm{decreases by $e(v,B) \geq
    k$}}\!\!\!\!-\; k|B| -
\underbrace{odd_k(A,B)\rule[-5mm]{0pt}{0pt}}_{\textrm{decreases by $\leq
    1$}}\!.$$
So by Theorem \ref{thm:5.2}\,(a), $A \cup \{v\},B$ is also a $k$-Tutte-pair
in $G$, contradicting the assumption that $A,B$ is a $k$-Tutte-pair with
$A$ as large as possible.

And for~(c), suppose that $d_G(t) \geq|A|+ k + odd_k(A,B) - 1$ for some
$t\in B$. This implies that $d_{G-A}(t)\geq k + odd_k(A,B) - 1$. Now
move~$t$ to $C$, and consider the change in each term in $\Theta_k(A,B)$:
$$k|A| +
 \underbrace{{\textstyle\sum\limits_{u \in B}}
   d_{G-A}(u)\rule[-5mm]{0pt}{0pt}}_{\substack{\textrm{decreases by} \\[1pt]
     d_{G-A}(t)\,\geq\, k+odd_k(A,B) - 1}}
\! -\!
\underbrace{k|B|\rule[-5mm]{0pt}{0pt}}_{\textrm{decreases by $k$}} - 
\underbrace{odd_k(A,B)\rule[-5mm]{0pt}{0pt}}_{\textrm{decreases by $\leq odd_k(A,B)$}}\!\!\!\!\!\!\!\!\!\!.$$
So by Theorem \ref{thm:5.2}\,(a), $A,B-\{t\}$ is also a $k$-Tutte-pair for
$G$, contradicting the minimality of $A,B$.\openeop

\bigskip
We introduce some further notation. Set $a \doteq |A|$, $b \doteq |B|$,
$c \doteq |C|=n-a-b$, $q\doteq odd_k(A,B)$, $r \doteq a + k+q-2$, and
$s\doteq n-\max\{0,b-k+1\}-\max\{0,q-1\}-1$. Using this notation,
\eqref{eqn12} can be written as
\begin{equation}\label{eqn12a}
\sum_{u\in B}d_G(u)\le kb+ab-ka+q-2.
\end{equation}
By Lemma \ref{lem:5.4}\,(a) we have $0\le a<\frac12n$. Since~$B$ is
disjoint from~$A$, we trivially have $0\le b\le n-a$. And since the number
of odd components of~$C$ is at most the number of elements of~$C$, we are
also guaranteed that $q\le n-a-b$. Finally, since for all vertices~$v$ we
have $d_G(v)\ge d_1\ge k$, we get from~\eqref{eqn12a} that
$q\ge\sum\limits_{u\in B}d_G(u)-kb-ab+ka+2\ge kb-kb-ab+ka+2=a(k-b)+2$,
hence $q\ge\max\{0,a(k-b)+2\}$. It follows that $a,b,q$ satisfy the
conditions in Theorem~\ref{thm:k-factor}\,(ii).

Next, by Lemma \ref{lem:5.4}\,(c) we have that
\begin{equation}\label{eqn13}
\text{for all $u\in B$:}\quad d_G(u) \le r.
\end{equation}

If $C\neq\varnothing$ (i.e., if $a+b<n$), let $m$ be the size of a largest
component of $\langle C\rangle$. Then, using Lemma \ref{lem:5.4}\,(b), for all
$v\in C$ we have
\begin{align*}
  d_G(v)&{}=e(v,A)+e(v,B)+e(v,C)\le|A|+\min\{k-1,|B|\}+m-1\\
  &{}=a+b-\max\{0,b-k+1\}+m-1.
\end{align*}
Clearly $m\le|C|=n-a-b$. If $q\ge1$, then $m\le n-a-b-(q-1)$, since~$C$ has
at least~$q$ components. Thus $m\le n-a-b-\max\{0,q-1\}$. Combining this
all gives
\begin{equation}\label{eqn14}
  \text{for all $v\in C$:}\quad d_G(v)\le
  n-\max\{0,b-k+1\}-\max\{0,q-1\}-1=s.
\end{equation}

Next notice that we cannot have $n-a=0$, because otherwise
$B=C=\varnothing$ and $odd_k(A,B)=0$, and \eqref{eqn12} becomes
$0\le-ka-2$, a contradiction. From~\eqref{eqn13} and~\eqref{eqn14} we see
that each of the $n-a>0$ vertices in $B \cup C$ has degree at most
$\max\{r,s\}$, and so $d_{n-a}\le\max\{r,s\}$.

If $r\le s$, then each of the $b$ vertices in~$B$ has degree at most~$r$,
and so $d_b\le r$. This also holds if $b=0$, since we set $d_0=0$,
and $r=a+k+q-2\ge0$ because $k\ge2$.

If $r>s$, then each of $n-a-b$ vertices in~$C$ has degree at most $s$
by~\eqref{eqn14}, and so $d_{n-a-b}\le s$. This also holds if $n-a-b=0$,
since we set $d_0=0$ and
\begin{align*}
  s&{}=n-\max\{0,b-k+1\}-\max\{0,q-1\}-1\\
&{}\ge\min\{n-1,n-q,(n-b)+(k-2),(n-q-b)+(k-1)\}\ge0,
\end{align*}
since $k\ge2$ and $q\le n-a-b$.

So we always have $r\le s$ and $d_b\le r$, or $r>s$ and $d_{n-a-b}\le s$,
but also $d_{n-a}\le\max\{r,s\}$, contradicting assumption (ii)\,$(*)$ in
Theorem~\ref{thm:k-factor}.\eop

\bigskip
How good is Theorem \ref{thm:k-factor}? We know it is not best monotone for
$k=2$. For example, the sequence $\pi = 4^4 6^3 10^4$ satisfies Theorem
\ref{thm:2-factor}, but not Theorem \ref{thm:k-factor} (it violates~$(*)$
when $a=4$, $b=5$ and $q=2$, with $r=6$ and $s=5$). And it is very unlikely
the theorem is best monotone for any $k \geq 3$. Nevertheless, Theorem
\ref{thm:k-factor} appears to be quite tight. In particular, we conjecture
for each $k \geq 2$ there exists a $\pi = (d_1 \leq \dots \leq d_n)$ such
that{

\qitem{$\bullet$} $(\pi,k)$ satisfies Theorem \ref{thm:k-factor}, and

\qitem{$\bullet$} there exists a degree sequence~$\pi'$, with
$\pi' \leq \pi$ and
$\sum\limits_{i=1}^n d_i' = \Bigl( \sum\limits_{i=1}^n d_i \Bigr) - 2$,
such that $\pi'$ is not forcibly $k$-factor graphical.

}Informally, for each $k \geq 2$, there exists a pair $(\pi,\pi')$ with
$\pi'$ `just below' $\pi$ such that Theorem \ref{thm:k-factor} detects that
$\pi$ is forcibly $k$-factor graphical, while $\pi'$ is not forcibly
$k$-factor graphical.

For example, let $n \equiv 2\pmod 4$ and $n \geq 6$, and consider the
sequences $\pi_n \doteq \bigl(\frac12 n \bigr)^{n/2 + 1} (n-1)^{n/2 - 1}$
and $\pi_n' \doteq \bigl(\frac12 n - 1 \bigr)^2 \bigl(\frac12 n \bigr)^{n/2
 - 1} (n-1)^{n/2-1}$. It is easy to verify that the unique realization of
$\pi_n'$ fails to have a $k$-factor, for $k = \frac14(n+2) \geq 2$. On the
other hand, we have programmed Theorem \ref{thm:k-factor}, and verified
that $\pi_n$ satisfies Theorem \ref{thm:k-factor} with $k=\frac14(n+2)$ for
all values of $n$ up to $n=2502$. We conjecture that
$(\pi_n, \frac14(n+2) )$ satisfies Theorem \ref{thm:k-factor} for all
$n \geq 6$ with $n \equiv 2\pmod4$.

\medskip
There is another sense in which Theorem \ref{thm:k-factor} seems quite
good. A graph~$G$ is \emph{$t$-tough} if $t\cdot\omega(G) \le |X|$, for
every $X\subseteq V(G)$ with $\omega(G-X)>1$, where $\omega(G-X)$ denotes
the number of components of $G-X$. In \cite{BBKSV2}, the authors give the
following best monotone condition for $\pi$ to be forcibly $t$-tough, for
$t \geq 1$.

\begin{thm}[\cite{BBKSV2}]\label{tough}
  \; Let $t\ge1$, and let $\pi=(d_1 \leq \dots \leq d_n)$ be graphical with
  $n > (t+1) \lceil t \rceil / t$. If
  $$d_{\lfloor i/t \rfloor} \leq i \; \Longrightarrow \;
  d_{n-i} \geq n-\lfloor i/t \rfloor,\qquad \text{for $t \le i <
    tn/(t+1)$},$$
  then $\pi$ is forcibly $t$-tough graphical.
\end{thm}

We also have the following classical result.

\begin{thm}[\cite{J4}]\label{J4-2}
  \; Let $k \geq 1$, and let $G$ be a graph on $n \geq k+1$ vertices with
  $kn$ even. If $G$ is $k$-tough, then $G$ has a $k$-factor.
\end{thm}

Based on checking many examples with our program, we conjecture that there
is a relation between Theorems \ref{tough} and \ref{thm:k-factor}, which
somewhat mirrors Theorem \ref{J4-2}.

\begin{conj}
\quad Let $\pi = (d_1 \leq \dots \leq d_n)$ be graphical, and let
$k \geq 2$ be an integer with $n> k+1$ and $kn$ even. If $\pi$ is forcibly
$k$-tough graphical by Theorem~\ref{tough}, then $\pi$ is forcibly
$k$-factor graphical by Theorem \ref{thm:k-factor}.
\end{conj}

\end{document}